\documentclass[12pt,reqno]{amsart}
  \usepackage[all,2cell]{xy}
  \usepackage{amsfonts}
  \usepackage{amsthm}
  \usepackage{amsmath}
  \usepackage{amssymb}
 \numberwithin{equation}{section}
 \usepackage{color}
\definecolor{darkgreen}{rgb}{0,0.45,0} 
\usepackage[colorlinks,citecolor=darkgreen,linkcolor=darkgreen]{hyperref}
  \usepackage{hyperref}
\usepackage{enumerate,xspace}

\UseAllTwocells
\CompileMatrices

\allowdisplaybreaks

\renewcommand{\epsilon}{\varepsilon}
\renewcommand{\phi}{\varphi}

\newcommand{\C}{\ensuremath{\mathcal C}\xspace}

\newcommand{\K}{\ensuremath{\mathcal K}\xspace}
\renewcommand{\L}{\ensuremath{\mathcal L}\xspace}
\newcommand{\M}{\ensuremath{\mathcal M}\xspace}

\newcommand{\Iso}{\textnormal{\bf Iso}\xspace}

\newcommand{\Ab}{\textnormal{\bf Ab}\xspace}

\newcommand{\Cat}{\textnormal{\bf Cat}\xspace}
\newcommand{\Mod}{\textnormal{\bf Mod}\xspace}

\newcommand{\MMod}{\textnormal{Mod(\M)}\xspace}

\newcommand{\op}{\ensuremath{^{\textnormal{op}}}}
\newcommand{\co}{\ensuremath{^{\textnormal{co}}}}
\newcommand{\coop}{\ensuremath{^{\textnormal{coop}}}}
\newcommand{\ob}{\textnormal{ob}}

\newcommand{\cl}{\colon}

\newcommand{\ls}{/\!/}

\DeclareMathOperator{\alg}{alg}
\DeclareMathOperator{\Lax}{Lax}
\DeclareMathOperator{\Mnd}{Mnd}

\newcommand{\LaxT}{\ensuremath{\textnormal{Lax-$T$-Alg}_c}\xspace}
\newcommand{\LaxG}{\ensuremath{\textnormal{Lax-$G$-Coalg}_\ell}\xspace}

\def\x{\times}

\def\1c#1{\stackrel{#1}{\to}}

  \newtheorem{proposition}{Proposition}[section]

  \newtheorem{theorem}[proposition]{Theorem}

  \theoremstyle{definition}
  \newtheorem{definition}[proposition]{Definition}
  \newtheorem{example}[proposition]{Example}

  \theoremstyle{remark}
  \newtheorem{remark}[proposition]{Remark}

  \newcounter{c}
  
  \newcommand{\etyk}[1]{\vspace{-7.4mm}$$\begin{equation}\Label{#1}
  \addtocounter{c}{1}}
  \renewcommand{\]}{\ifnum \value{c}=1 $$\else \end{equation}\fi}
\begin{document}

 \title{Morita contexts as lax functors}

\author{Stephen Lack}
\thanks{This research was supported by the Australian Research
Council in the form of a Future Fellowship as well as a Discovery 
Project Grant, project number DP1094883.}
\address{Department of Mathematics\\ Macquarie University}
\email{steve.lack@mq.edu.au}

\date{\today}

\begin{abstract}
Monads are well known to be equivalent to lax functors out of the
terminal category. Morita contexts are here shown to be lax functors out of the chaotic category with two objects. This allows various aspects in the theory of Morita contexts to be seen as special cases of general results about lax functors. The account we give of this could serve as an introduction to lax functors for those familiar with the theory of monads. We also prove some very general results along these lines relative to a given 2-comonad, with the 
classical case of ordinary monad theory amounting to the case of the identity comonad on \Cat. 
\end{abstract}
  
\maketitle


\section{Introduction}

Some of the fundamental aspects of the theory of monads are:
\begin{enumerate}
\item the fact that every adjunction generates a monad, 
\item the construction of the Eilenberg-Moore category of a monad, and the fact that this is part of an adjunction generating the given monad,
\item the universal property of this construction,
\item Beck's theorem characterizing which adjunctions arise in this way.
\end{enumerate}
All of these can be formulated and proved in a suitable 2-category; 
the resulting theory is known as the {\em formal theory of monads} 
\cite{ftm}.

The recent paper \cite{BrzezinskiMarquezVercruysse} formulated
and proved analogous results involving {\em Morita contexts} rather 
than monads. A Morita context, in the sense of \cite{BrzezinskiMarquezVercruysse}, involves a monad $t$ on a
category $A_0$, a monad $s$ on a category $A_1$, functors 
$f\cl A_0\to A_1$ and $g\cl A_1\to A_0$, and various further natural 
transformations subject to compatibility conditions. The 
corresponding notion of adjunction consists of functors 
$u_0\cl B\to A_0$ and $u_1\cl B\to A_1$ with the same domain, each 
possessing a left adjoint. These were studied in 
\cite{BohmMenini-pretorsors} under the name {\em double
adjunction}. The word ``double'' refers to the fact that there are
two adjunctions, not to any connection with double categories. 

The main purpose of this note is to point out a common 
generalization of monads and Morita contexts, in which the 
various results listed above are all known. 

B\'enabou observed in  \cite{bicategories} that to give a monad 
is equivalently to give a {\em lax functor} from the terminal 
2-category 1 to \Cat. Building on this, Street \cite{Street-twoconstructions} studied lax 
functors from $X$ to \Cat for an arbitrary small category $X$,
and managed to formulate and prove results which, when 
specialized to the case $X=1$, gave each of the results about
monads mentioned above. 

After recalling the notion of lax functor, we shall quickly focus
on the case where $X$ is the category \Iso with two objects $0$ 
and $1$, and exactly one arrow in each of the four hom-sets. This is
often called the ``chaotic'' or ``indiscrete'' category on two objects; it
is also called the ``free-living isomorphism'', since 
a functor from \Iso to a category $D$ is precisely an isomorphism
in $D$. The key observation of this paper is that a {\em lax functor
from \Iso to \Cat is precisely a Morita context.} 
After explaining this, we 
then recall aspects  of \cite{Street-twoconstructions} and explain how they can be used to obtain the main theoretical results 
in \cite{BrzezinskiMarquezVercruysse}.

There is also another bicategorical treatment of Morita contexts, introduced in \cite{ElKaoutit}, and going under the name of {\em wide Morita context}. These wide Morita contexts are the special case of the Morita contexts studied here and in \cite{BrzezinskiMarquezVercruysse} in which the two monads are both trivial (identity monads). Corresponding to our observation that Morita contexts in \M are the same as lax functors from \Iso to \M is the observation made in \cite{Pecsi} that {\em wide} Morita contexts in \M are the same as {\em normal} lax functors from \Iso to \M, where a lax functor is said to be normal when it strictly preserves the identities.

Thus wide Morita contexts in a bicategory \M are a special case of Morita contexts in \M; on the other hand, there is also a way to see Morita contexts as a special case of wide Morita contexts, as we now explain. For a bicategory \M whose hom-categories have reflexive coequalizers, preserved by composition on either side, there is a bicategory \MMod having monads in \M as its objects, and ``2-sided actions'' as 1-cells: when \M is the one-object bicategory \M corresponding to the monoidal category \Ab of abelian groups, the corresponding \MMod is the bicategory \Mod of rings, bimodules, and homomorphisms of bimodules. There is a natural bijection between lax functors with codomain \M and normal lax functors with codomain \MMod; taking the case where the domain is \Iso, we obtain a bijection between Morita contexts in \M and wide Morita contexts in \MMod. A classical Morita context between rings is a Morita context in the one-object bicategory corresponding to \Ab; or, equivalently, a wide Morita context in \Mod.

There is a ``formal theory of lax functors with domain $X$'', 
dealing with lax functors from $X$ to a general 2-category \M,
and containing the formal theory of monads as the special case
$X=1$, but we shall restrict ourselves to the case $\Cat$. 
On the other hand in the final two sections we describe a still more 
general setting using the language of 2-dimensional monad theory \cite{BKP}, which includes the case of lax functors from an arbitrary small 2-category $X$ to a complete 2-category \M. 
Readers who enjoy the adrenaline rush of Extreme 2-Category Theory are welcome to skip straight to Section~\ref{sec:generalizations}.

For an object $x$ of a category or 2-category we shall often
write $x$ for the corresponding identity 1-cell; similarly 
identity 2-cells in 2-categories will sometimes be denoted by 
the name of the corresponding 1-cell. Alternatively, we may 
just write $1$ for an identity if the domain/codomain is clear
from the context. 

We shall often consider a category, such as $X$, as a 
2-category with no non-identity 2-cells.


\section{Monads and lax functors}

A lax functor from a category $X$ to a 2-category \K  consists of 
the following assignments. For each object $x\in X$ there is a 
specified object $A_x\in\K$, for each morphism $\xi\cl x\to y$
in $X$ there is a morphism (1-cell) $t_\xi\cl A_x\to A_y$ in \K,
for each composable pair of morphisms $\xi\cl x\to y$ 
and $\zeta\cl y\to z$ there is a 2-cell 
$\mu_{\xi,\zeta}\cl t_\zeta t_\xi\to t_{\zeta\xi}$ in \K, 
and for each object $x\in X$ there is a 2-cell
$\eta_x\cl 1_{A_x}\to t_x = t_{1_x}$ in \K. We shall sometimes omit
the subscripts on $\mu$ and $\eta$.
These are required to satisfy the
associativity condition stating that the diagram
$$\xymatrix{
t_\tau t_\zeta t_\xi \ar[r]^{\mu t_\xi} 
\ar[d]_{t_\tau \mu} & 
t_{\tau\zeta} t_\xi \ar[d]^{\mu} \\
t_\tau t_{\zeta\xi} \ar[r]_{\mu} & t_{\tau\zeta\xi} 
}$$
commutes for all composable triples $\xi,\zeta,\tau$.
They are also required to satisfy the unit conditions asserting
the commutativity of the triangles
$$\xymatrix{
t_\xi 1 \ar[r]^{t_\xi\eta} \ar@{=}[dr] & 
t_\xi t_1 \ar[d]^{\mu} & 
1 t_\xi \ar[r]^{\eta t_\xi} \ar@{=}[dr] & 
t_1 t_\xi \ar[d]^{\mu}  \\
& t_\xi && t_\xi
}$$
for all arrows $\xi$.

The slightly idiosyncratic choice of notation is of course 
designed to emphasize that if $X=1$ then this is just 
a monad. 

Now consider the case  $X=\Iso$. We write $\alpha\cl x\to y$
and $\beta\cl y\to x$ for the two non-identity arrows of $X$.
There are exactly 8 composable pairs and exactly 16 composable
triples. A lax functor $\Iso\to\Cat$ involves
categories $A_x$ and $A_y$, a monad $(t_x,\mu,\eta)$
on $A_x$ and a monad $(t_y,\mu,\eta)$ on $A_y$. There 
are also functors $t_\alpha\cl A_x\to A_y$ and $t_\beta\cl A_y\to A_x$
which we shall call $f$ and $g$ respectively. There are 
6 further 2-cells, corresponding to the 6 remaining composable
pairs in \Iso; explicitly, they are $\rho_f=\mu_{\alpha,y}\cl t_y f\to f$,
$\lambda_f=\mu_{x,\alpha}\cl f t_x\to f$, 
$\rho_g=\mu_{x,\beta}\cl g t_y\to g$, 
$\lambda_g= \mu_{\beta,y}\cl  t_x g\to g$,
$\epsilon_x=\mu_{\alpha,\beta}\cl gf\to t_x$,
and $\epsilon_y=\mu_{\beta,\alpha}\cl fg\to t_y$. There are 14 
further associativity conditions and 4 further unit conditions,
corresponding precisely to the conditions in the definition 
of Morita context:  see \cite{BrzezinskiMarquezVercruysse}.

\section{Algebras and lax natural transformations}

Let $A\cl X\to\Cat$ be a lax functor, involving objects $A_x$,
morphisms $t_\xi$, and 2-cells $\mu_{\xi,\zeta}$ and 
$\eta_x$ as above; and let $B\cl X\to\Cat$ be another lax
functor, involving $B_x$, $s_\xi$, $\mu_{\xi,\zeta}$, and $\eta_x$.

A lax natural transformation $u$ from $B$ to $A$ consists of 
a morphism $a_x\cl B_x\to A_x$ for each $x\in X$ and a 2-cell
$\alpha_\xi$ as in 
$$\xymatrix{
B_x \ar[r]^{a_x}_{~}="1" \ar[d]_{s_\xi} & 
A_x \ar[d]^{t_\xi} \\
B_x \ar[r]_{a_y}^{~}="2" & A_d
\ar@{=>}"1";"2"^{\alpha_\xi} }$$
for each morphism $\xi\cl x\to y$ in $X$. These are required
to be compatible with the $\mu$ and with the $\eta$. The 
precise conditions can be found for example in 
\cite[Section~1]{Street-twoconstructions}, where the name
{\em left lax transformation} is used. 

In the case where $X=1$, so that a lax functor $X\to\Cat$ is 
just a monad, this reduces to a functor $a\cl B\to A$ equipped with 
a 2-cell $\alpha\cl ta\to as$ satisfying the two conditions for 
$(a,\alpha)$ to be a monad morphism from $(B,s)$ to $(A,t)$, 
in the sense of \cite{ftm}.

Specializing further, if $B$ is strict, so that $s=1$, this reduces
to a functor $a\cl B\to A$ equipped with an action $\alpha\cl ta\to a$;
and if finally $B$ is the terminal category, this is just an object
$a\in A$ equipped with a $t$-algebra structure $\alpha\cl ta\to a$.

Reverting to the general case, 
there is also a notion of {\em oplax natural transformation},
called {\em right lax transformation} in \cite{Street-twoconstructions},
in which the direction of the $\alpha_\xi$ are reversed. In 
the case $X=1$ these are the {\em monad opfunctors} of 
\cite{ftm}. 

There are also {\em modifications} between lax natural transformations, which capture the monad 2-cells of \cite{ftm}. A modification
between lax natural transformations $(a,\alpha)$ and $(b,\beta)$ as above consists of a natural transformation
$\phi_x\cl a_x\to b_x$ for each $x\in X$ suitably compatible with 
$\alpha$ and $\beta$: see \cite{Street-twoconstructions}. 
If $X=1$, and $B$ is the strict functor with value the terminal 
category, so that $(a,\alpha)$ and $(b,\beta)$ are just $t$-algebras,
then such a $\phi$ is precisely a morphism of $t$-algebras.

There is a 2-category $\Lax(X,\Cat)_\ell$ of lax functors, lax natural 
transmorphisms, and modifications which reduces, in the case
$X=1$, to Street's 2-category $\Mnd(\Cat)$ of monads in \Cat.

\section{Adjunctions}

Let $B\cl X\to\Cat$ be a lax functor, and suppose that for each
object $x\in X$ we are given a category $A_x$, a functor 
$u_x\cl B_x\to A_x$ and a left adjoint $f_x\dashv u_x$ with unit
$\sigma_x$ and counit $\epsilon_x$. 

For each $\xi\cl x\to y$ in $X$, let $s_\xi\cl A_x\to A_y$
be given by the composite 
$$\xymatrix{
A_x \ar[r]^{f_x} & B_x \ar[r]^{t_\xi} & B_y \ar[r]^{u_y} & A_y.}$$
For $\xi$ as above and $\zeta\cl y\to z$, define
$\mu_{\xi,\zeta}$ for $A$ to be the composite
$$\xymatrix{
&&& A_y \ar[dr]_{f_y} \ar[ddrrr]^{s_\zeta} \\
&& B_y \ar[ur]_{u_y} \ar[rr]_{1}="2"^{~}="1" \ar@{=>}[ur];"1"^{\epsilon_y} &&
B_y \ar[dr]_{t_\zeta} \\
A_x \ar[r]^{f_x} \ar[uurrr]^{s_\xi} \ar@/_2pc/[rrrrrr]_{s_{\zeta\xi}} & 
B_x \ar[ur]_{t_\xi} \ar[rrrr]_{t_{\zeta\xi}}^{~}="3" &&&& 
B_z \ar[r]^{u_z} &  A_z 
\ar@{=>}"2";"3"^{\mu_{\xi,\zeta}} 
}$$
and define the $\eta$ for $A$ to be the composite
$$\xymatrix{
A_x \ar[r]_{f_x} \ar@/^3pc/[rrr]^{1}_{~}="1" \ar@/_2pc/[rrr]_{s_x} &
B_x \ar@/^1pc/[r]^{1}="2"_{~}="3" \ar@/_1pc/[r]_{t_x}^{~}="4" & 
B_x \ar[r]_{u_x} & A_x 
\ar@{=>}"1";"2"^{\sigma_x}
\ar@{=>}"3";"4"^{\eta_x}
}$$
\begin{proposition}[Street]
This defines a lax functor $A$;  the $f_x$ become 
the components of an oplax natural transformation and the
$u_x$ become the components of a lax natural transformation.  
\end{proposition}

The most important case is the following.

\begin{definition}
  An $X$-adjunction consists of a strict functor $B:X\to\Cat$
equipped with an adjunction $f_x\dashv u_x\cl B_x\to A_x$ for 
each object $x\in X$.
\end{definition}

If $X=1$, then $B$ reduces to a single category,
and we are given a functor $u\cl B\to A$ with a left adjoint $f\dashv u$, so that a $1$-adjunction is just an adjunction.
The induced lax functor $A\cl 1\to\Cat$ is the monad induced by 
the adjunction. 

In the case $X=\Iso$, then $B$ amounts to a pair of categories
with an isomorphism between them, but up to isomorphism this
is again just a single category. But our notion of adjunction is 
now a pair of categories $A_x$ and $A_y$, and a pair of 
functors $u_x\cl B\to A_x$ and $u_y\cl B\to A_y$ with adjunctions
$f_x\dashv u_x$ and $f_y\dashv u_y$; thus an $\Iso$-adjunction 
amounts to a 
double adjunction in the sense of \cite{BohmMenini-pretorsors}.
The induced lax functor $\Iso\to\Cat$ is precisely the 
induced Morita context described in \cite[Section~3.1]{BrzezinskiMarquezVercruysse}.


\section{Eilenberg-Moore construction}

We have seen that a monad is the same thing as a lax functor
$1\to\Cat$, and that an algebra for the monad is then the same
thing as a lax natural transformation from the strict functor
$1\to\Cat$ whose image is the terminal category. This strict
functor is of course the unique representable functor
$1\to\Cat$. This 
suggests the importance of lax natural transformations from 
strict functors, especially representable ones,  to lax functors. 

Let $[X,\Cat]$ denote the 2-category of strict functors 
from $X$ to \Cat, strict natural transformations between them,
and modifications. Recall that we write $\Lax(X,\Cat)_\ell$ 
for the 2-category of lax functors from $X$ to \Cat, lax natural
transformations, and modifications.

\begin{proposition}[Street]
The inclusion $[X,\Cat]\to\Lax(X,\Cat)_\ell$ has a right adjoint.   
\end{proposition}

We shall write $\alg(A)\cl X\to\Cat$ for the value of the right
adjoint at a lax functor $A\cl X\to\Cat$. The universal property 
of the adjunction says in part that for any strict functor
$B\cl X\to\Cat$, there is a bijection between lax natural
transformations $B\to A$ and strict natural transformations
$B\to\alg(A)$. In particular, this should be the case if
$B$ is a representable functor $X(x,-)\cl X\to\Cat$. By the Yoneda
lemma, strict natural transformations $X(x,-)\to\alg(A)$ 
are in bijection with objects of $\alg(A)_x$, and so we deduce
that the objects of $\alg(A)_x$ should be the lax natural 
transformations from $X(x,-)$ to $A$. Similarly the morphisms
of $\alg(A)_x$ are the modifications between lax natural transformations from $X(x,-)$ to $A$.

This reasoning tells us what $\alg(A)$ must be; the fact that this
actually works is proved in \cite[Theorem~4]{Street-twoconstructions}.

 If $X=1$, so that a lax functor $A\cl 1\to\Cat$ is just a monad,
then $\alg(A)\cl 1\to\Cat$ picks out the Eilenberg-Moore category 
of the monad. 
This $\alg(A)$ is Street's ``second construction'' on lax 
functors; the first corresponds to the Kleisli construction
for monads, and provides a {\em left} adjoint to the inclusion of
$[X,\Cat]$ in the 2-category $\Lax(X,\Cat)_{c}$ of lax functors,
{\em oplax} natural transformations, and modifications.

\begin{proposition}[Street]
   Let $A\cl X\to\Cat$ be a lax functor, and $\alg(A)\cl X\to\Cat$ the 
strict functor defined above. For each $x\in X$ there is a 
functor $u_x\cl \alg(A)_x\to A_x$ sending a lax natural transformation
$X(x,-)\to A$ to the value at $1\cl x\to x$ of its $x$-component
$X(x,x)\to A_x$, and this functor $u_x$ has a left adjoint $f_x$. 
The lax functor $X\to\Cat$ induced by this adjunction is $A$.
\end{proposition}

For a Morita context, seen as a lax functor $A\cl \Iso\to\Cat$, the 
corresponding strict functor $\alg(A)\cl \Iso\to\Cat$ will
give a pair of isomorphic categories, but up to isomorphism
we may take this just to be a single category. This is exactly 
the {\em Eilenberg-Moore category associated to the Morita
context}, as described in \cite[Section~3.2]{BrzezinskiMarquezVercruysse}.

\section{Interlude: parametrized internal coequalizers}
\label{sec:pic}

The Beck condition characterizing monadic functors involves
existence and preservation of certain coequalizers. We shall
describe a generalization of Beck's theorem in the next section,
but first we need to develop a notion of coequalizer suitable
for our purposes. 

Let \K be a 2-category, and $K$ an object of \K. We write 
$\K\ls K$ for the evident 2-category whose objects are arrows
in \K with codomain $K$, and arrows are triangles in \K which 
commute up to a specified 2-cell. We use the following notation.
An object of $\K\ls K$ is a morphism $A\cl\partial_0A\to K$ in \K.
A morphism $A\to B$ consists of a morphism 
$\partial_0f\cl\partial_0B\to\partial_0A$ in \K equipped with a 
2-cell $f$ from $A$ to the composite of $B$ and $\partial_0 f$. A
2-cell from $f$ to $g$ is a 2-cell $\partial_0 f\to \partial_0 g$
which when pasted with $g$ gives $f$.


There is an evident forgetful 2-functor $\partial_0\cl \K\ls K\to\K$.
The image under $\partial_0$ of an object or morphism will sometimes be called its {\em head}. 

\begin{example}
  If $\K=\Cat$, and $K\in\Cat$ is a category, then we may 
identify $K$ with the full subcategory of $\K\ls K$ consisting 
of all objects $(A,a)$ with $A$ the terminal category 1. In 
other words, $K$ is the fibre of $\partial_0\colon\Cat\ls K\to\Cat$ over $1$.
\end{example}

\begin{definition}
We say that the object $K\in\K$ has {\em parametrized internal
coequalizers} if the category $\K\ls K$ has coequalizers preserved
by $\partial_0$.   
\end{definition}

We shall often want to restrict to the case of parametrized 
internal coequalizers with specified head; in particular, 
the example above shows that a category $K$ has ordinary
coequalizers if and only if it has parametrized internal coequalizers
in \Cat with head the terminal coequalizer diagram 
$$\xymatrix{
1 \ar@<1ex>[r] \ar@<-1ex>[r] & 1 \ar[r] & 1. }$$
 
Since our parametrized internal coequalizers in $K$ are just ordinary 
coequalizers in $\K\ls K$, we can define a parametrized internal
coequalizer to be split when it is split as a coequalizer in 
$\K\ls K$. Once again, we may choose to specify the splitting 
in the head.

When it comes to preservation, both the 2-category \K and 
the object $K$ might be varied. The key observation is that 
for a 2-functor $F\cl \K\to\L$ and a morphism $f\cl FK\to L$ there
is an induced 2-functor $F\ls f\cl \K\ls K\to \L\ls L$ sending 
an object $A\cl \partial_0 A\to K$ to the object 
$$\xymatrix{
F(\partial_0A) \ar[r]^{FA} & FK \ar[r]^{f} &  L }$$
of $\L\ls L$. We shall say that $F\ls f$ {\em preserves} a parametrized
internal coequalizer in $\K\ls K$ if $F\ls f$ sends it to a parametrized 
internal coequalizer in $\L\ls L$. This includes in particular
the requirement that $F$ preserve the coequalizer appearing 
in the head.

\section{Beck theorem}

Consider an $X$-adjunction consisting of a strict functor
$B\cl X\to\Cat$ equipped with functors
$u_x\cl B_x\to A_x$ for each $x\in X$,  with left adjoints 
$f_x\dashv u_x$. 
We have seen that these adjunctions generate a lax functor
$A\cl X\to\Cat$, and that the $u_x$ become the components of 
a lax natural transformation $u\cl B\to A$. Thus by the universal
property of the Eilenberg-Moore object (second construction)
they correspond to a unique strict natural transformation 
$k\cl B\to\alg(A)$. 

In the case $X=1$, this is the canonical comparison functor
into the Eilenberg-Moore category; in the case $X=\Iso$ it
is the comparison functor of 
\cite[Section~3.5]{BrzezinskiMarquezVercruysse}.

We shall say that the $X$-adjunction is {\em monadic} if the 
components $k_c\cl B_c\to\alg(A)_c$ of $k$ are equivalences
of categories. This is not enough to make $k$ into an equivalence
in the 2-category $[X,\Cat]$, since the equivalence inverses
$\alg(A)_c\to B_c$ need not be strictly natural, but they will be
pseudonatural, and in fact the condition is equivalent to saying 
that $k$ is an equivalence in the 2-category of 2-functors
from $X$ to \Cat, pseudonatural transformations, and modifications.

This corresponds to the ordinary (up-to-equivalence) notion of 
monadicity in the case $X=1$, and to what was called 
``moritability'' in the case $X=\Iso$. 
\cite{BrzezinskiMarquezVercruysse}

Similarly, we may define the $X$-adjunction to be {\em strictly
monadic} if the $k_c$ are isomorphisms of categories. This is
of course equivalent to $k$ being an isomorphism in $[X,\Cat]$. 

It is this strict notion of monadicity that was used in 
\cite{Street-twoconstructions}. The Beck-style theorem of
\cite[Theorems~10~and~11]{Street-twoconstructions} involved
a colimit-notion called a ``universal reflection of a centipede'',
which we shall now describe using the notion of parametrized 
internal coequalizer given in the previous section. The presentation
given here looks quite different to that of \cite{Street-twoconstructions} but the equivalence of the two presentations is an extended but straightforward exercise using the Yoneda lemma.

Write $\ob X$ for the category with the same objects as $X$,
but with only identity morphisms. The evident forgetful 
2-functor $U\cl [X,\Cat]\to[\ob X,\Cat]$ given by composition with
the inclusion $\ob X\to X$ has both left and right adjoints, given
by left and right Kan extension, and in fact $U$ is both monadic
and comonadic. In particular, we shall write $F$ for the left adjoint
of $U$, and $\epsilon\cl FU\to 1$ for the counit and $\eta\cl 1\to UF$
for the unit. A straightforward calculation shows that, for an 
arbitrary strict functor $C\cl X\to\Cat$, the functor $FUC$ is 
given by 
$$(FUC)_x = \sum_{y\to x} C_y$$
with $\epsilon\cl FUC\to C$ induced by the maps $C_\xi\cl C_y\to C_x$
for $\xi\cl y\to x$. One aspect of the monadicity of $U$ is that 
the diagram 
$$\xymatrix{
FUFUC \ar@<1ex>[r]^{\epsilon FU} \ar@<-1ex>[r]_{FU\epsilon} & 
FUC \ar[r]^{\epsilon} & C  }$$
is always a coequalizer, and this coequalizer is $U$-split,
in the sense that it is sent by $U$ to a split coequalizer,
given in the following diagram.
\begin{equation}
  \label{eq:canonical-splitting}
  \xymatrix{
UFUFUC \ar@<1ex>[r]^{U\epsilon FU} \ar@<-1ex>[r]_{UFU\epsilon} & 
UFUC \ar[r]^{U\epsilon} \ar@/^2pc/[l]^{\eta UFU} & 
UC \ar@/^1pc/[l]^{\eta U}  }
\end{equation}

For a functor $B\cl X\to\Cat$ and object $x\in X$, an {\em 
$x$-centred centipede in $B$}, 
in the sense of \cite{Street-twoconstructions},
is now a pair of parallel arrows in $[X,\Cat]\ls B$ with head given 
by the parallel arrows in  
$$\xymatrix{
FUFUX(x,-) \ar@<1ex>[r]^{\epsilon FU} \ar@<-1ex>[r]_{FU\epsilon} & 
FUX(x,-) \ar[r]^{\epsilon} & X(x,-)  }$$
while a {\em universal reflection} for the centipede is precisely
an internal parametrized coequalizer; of course its head must 
be given as in the diagram, since we require the coequalizer
to be preserved by $\partial_0\cl [X,\Cat]\ls B\to[X,\Cat]$. An $x$-centred
centipede as above is {\em split by the family $u_x$}, in the 
sense of \cite{Street-twoconstructions}, just when it is sent 
by $U\ls u$ to a split coequalizer in $[\ob X,\Cat]\ls A$ with 
head given by 
the canonical splitting
appearing in \eqref{eq:canonical-splitting}. We shall say then
that the centipede is {\em $U\ls u$-split}.

We may now say that {\em $U\ls u$ creates ($U\ls u$)-split
parametrized internal coequalizers with canonical head}, if 
for each parallel pair of morphisms in $[X,\Cat]\ls B$ 
with head as in \eqref{eq:canonical-splitting} and which 
is $U\ls u$-split, there is a unique parametrized internal 
coequalizer in 
$[X,\Cat]\ls B$ which is sent by $U\ls u$ to the specified (split)
parametrized internal coequalizer in $[\ob X,\Cat]\ls A$.

\begin{theorem}[Street]
For an $X$-adjunction $f\dashv u\cl UB\to A$, the induced 
comparison map $k\cl B\to\alg(A)$ is invertible if and only if
$U\ls u$ creates ($U\ls u$)-split parametrized internal 
coequalizers with canonical head.
\end{theorem}

It is straightforward to modify this condition to deal with 
the ``up-to-equivalence'' notion of monadicity. One simply 
asks that $[X,\Cat]\ls B$ have, and that $U\ls u$ preserve and
reflect, parametrized
internal coequalizers with canonical head for all parallel pairs
whose image under $U\ls u$ has a split parametrized internal
coequalizer with canonical  head.

We now investigate what this condition says in the case 
$X=\Iso$ relevant to Morita contexts; as usual, we identify
strict functors $\Iso\to\Cat$ with single categories. 

An $X$-adjunction then consists of a pair of right adjoints
$u_0\cl B\to A_0$ and $u_1\cl B\to A_1$. 

A centipede consists of the following data:
objects $b_x\in B$ for each $x\in X$, 
objects $b_{xy}\in B$ for each $x,y\in X$,
morphisms $\beta_{xy}\cl b_{xy}\to b_x$ and $\beta'_{xy}\cl b_{xy}\to b_y$
for each $x,y\in X$. A universal reflection of the centipede is a
colimit of the diagram 
$$\xymatrix @R4pc @C3pc {
b_{00} \ar@<1ex>[d]^(0.3){\beta'_{00}} \ar@<-1ex>[d]_(0.3){\beta_{00}} &
b_{01} \ar@<1ex>[drrr]^(0.2){\beta'_{01}} \ar@<-1ex>[dl]_(0.2){\beta_{01}} &&
b_{10} \ar@<1ex>[dr]^(0.2){\beta'_{10}} \ar@<-1ex>[dlll]_(0.2){\beta_{10}} &
b_{11} \ar@<1ex>[d]^(0.3){\beta'_{11}} \ar@<-1ex>[d]_(0.3){\beta_{11}} \\
b_0 &&&& b_1 }$$
(compare diagram (16) in \cite{BrzezinskiMarquezVercruysse}).
The general theorem stated above asks for such universal 
reflections whenever the centipede is split by $U\ls u$, but on
inspecting the proof one finds, just as in the case of the classical
Beck theorem, that only certain specific colimits involving free 
algebras are required. It is these specific colimits which are
used in \cite[Theorem~3.12]{BrzezinskiMarquezVercruysse}.
 

\section{Generalizations}
\label{sec:generalizations}

As was mentioned in the previous section, the forgetful
2-functor $U\cl [X,\Cat]\to[\ob X,\Cat]$ is both monadic and
comonadic. In this section we focus on the comonadicity 
of $U$ along with the fact that the comonad in question preserves
certain limits (in the case of the previous section, it has an adjoint). 

We therefore consider a 2-category \K equipped with a 2-comonad
$G=(G,d,e)$, and we write $\K^G$ for the Eilenberg-Moore 
2-category of $G$: the 2-category of strict $G$-algebras,
strict morphisms of $G$-algebras, and algebra 2-cells.
Various weakenings of these notions have been studied in detail
in the monad case, but much less so for comonads, so we shall
recall in full the required definitions; this also serves to fix 
our conventions for the directions of 2-cells. We discuss 
in the next section how to dualize these results, and so 
obtain results about 2-monads. 

\begin{example}
The setting of \cite{Street-twoconstructions} corresponds to the 
case $\K=[\ob X,\Cat]$, with $G$ the comonad induced by 
right Kan extension. The classical case of ordinary monads 
corresponds to the case where $\K=\Cat$ and $G$ is the 
identity comonad. The setting of \cite{BrzezinskiMarquezVercruysse}
corresponds to the comonad $G$ on $\Cat\x\Cat$ induced by
the diagonal $\Delta\colon\Cat\to\Cat\x\Cat$ and its right
adjoint; explicitly $G(A,B)=(A\x B,A\x B)$.
\end{example}

A {\em lax $G$-coalgebra}
consists of an object $A$ equipped with a morphism $a\cl A\to GA$
and 2-cells
$$\xymatrix{
A \ar[r]^{a}_{~}="1" \ar[d]_{a} & GA \ar[d]^{Ga} &
A \ar[d]_{a} \ar[dr]^1_(0.4){~}="4" &  \\
GA \ar[r]_{dA}^{~}="2" & G^2A & GA \ar[r]_{eA}^(0.3){~}="3" & A 
\ar@{=>}"1";"2"^{\alpha}
\ar@{=>}"4";"3"^{\alpha_0} }$$
subject to the following three coherence conditions. 
\begin{align*}
  \vcenter{
\xymatrix{
& GA \ar[r]^{Ga} \ar[dr]_{Ga} \ar@{=>}[dd]^{\alpha} & 
G^2A \ar[dr]^{G^2a} \ar@{=>}[d]^{G\alpha} \\
A \ar[ur]^{a} \ar[dr]_{a} && G^2A \ar[r]^{GdA} & G^3A \\
& GA \ar[r]_{dA} \ar[ur]^{dA} & G^2A \ar[ur]_{dGA} 
}}
&=
  \vcenter{
\xymatrix{
& GA \ar[r]^{Ga} \ar@{=>}[d]^{\alpha} & G^2A \ar[dr]^{G^2a} \\
A \ar[ur]^{a} \ar[dr]_{a} \ar[r]^{a} & 
GA \ar[ur]_{dA} \ar[dr]^{Ga} \ar@{=>}[d]^{\alpha} && G^3A \\
& GA \ar[r]_{dA} & G^2A \ar[ur]_{dGA} 
}}
\\
\vcenter{
\xymatrix{
& GA \ar[dr]_{Ga} \ar@/^1pc/[drr]^{1}_{~}="1" \ar@{=>}"1";[dr]^{G\alpha_0}  \ar@{=>}[dd]^{\alpha} \\
A \ar[ur]^{a} \ar[dr]_{a} && G^2A \ar[r]_{GeA} & GA \\
& GA \ar[ur]^{dA} \ar@/_1pc/[urr]_{1} }}
&=
\vcenter{
\xymatrix{
& GA  \ar@/^1pc/[drr]^{1}_{~}="1" \\
A \ar[ur]^{a} \ar[dr]_{a} &&& GA \\
& GA  \ar@/_1pc/[urr]_{1} }}
\\
\vcenter{
\xymatrix{
& & A \ar[dr]^{a} \\
A \ar[r]^{a} \ar[dr]_{a} \ar@/^1pc/[urr]^{1}_{~}="1" \ar@{=>}"1";[r]^{\alpha_0} & GA \ar[ur]_{eA} \ar[dr]^{Ga} \ar@{=>}[d]^{\alpha} && GA \\
& GA \ar[r]_{dA} & G^2A \ar[ur]_{eGA} 
}}
&=
\vcenter{
\xymatrix{
& & A \ar[dr]^{a} \\
A  \ar[dr]_{a} \ar@/^1pc/[urr]^{1}_{~}="1"  & && GA \\
& GA \ar[r]_{dA} \ar@/^1pc/[urr]^1 & G^2A \ar[ur]_{eGA} 
}}
\end{align*}

A {\em lax
$G$-morphism} from $(A,a,\alpha,\alpha_0)$ to another
lax $G$-algebra $(B,b,\beta,\beta_0)$ consists of a morphism
$f\cl A\to B$ equipped with a 2-cell 
$$\xymatrix{
A \ar[r]^{f}_{~}="1" \ar[d]_{a} & B \ar[d]^{b} \\
GA \ar[r]_{Gf}^{~}="2" & GB 
\ar@{=>}"1";"2"^{\overline{f}} }$$
subject to two coherence conditions
\begin{align*}
  \vcenter{
\xymatrix{
& B \ar[r]^{b} \ar[dr]_{b} \ar@{=>}[dd]^{\overline{f}} & 
GB \ar[dr]^{Gb} \ar@{=>}[d]^{\beta} \\
A \ar[ur]^{f} \ar[dr]_{a} && GB \ar[r]_{dB} & G^2B \\
& GA \ar[r]_{dA} \ar[ur]_{Gf} & G^2A \ar[ur]_{G^2f} 
}}
&=
  \vcenter{
\xymatrix{
& B \ar[r]^{b} \ar@{=>}[d]^{\overline{f}} & GB \ar[dr]^{Gb} \ar@{=>}[dd]^{G\overline{f}} \\
A \ar[ur]^{f} \ar[dr]_{a} \ar[r]^{a} & GA \ar[ur]_{Gf} \ar[dr]^{Ga} \ar@{=>}[d]^{\alpha} && G^2B \\
& GA \ar[r]_{dA} & G^2A \ar[ur]_{G^2f} 
}}
\\
  \vcenter{
\xymatrix{
& B \ar@/^1pc/[drr]^{1}_{~}="1" \ar@{=>}"1";[dr]^{\beta_0} \ar[dr]_{b} \ar@{=>}[dd]^{\overline{f}} & \\
A \ar[ur]^{f} \ar[dr]_{a} && GB \ar[r]_{eB} & B \\
& GA \ar[r]_{eA} \ar[ur]_{Gf} & A \ar[ur]_{f} 
}}
&=
  \vcenter{
\xymatrix{
& B \ar@/^1pc/[drr]^{1}_{~}="1"  & \\
A \ar[ur]^{f} \ar[dr]_{a} \ar@/^1pc/[drr]^{1}_{~}="2" \ar@{=>}"2";[dr]^{\alpha_0} && & B \\
& GA \ar[r]_{eA}  & A \ar[ur]_{f} 
}}
\end{align*}

There are also
{\em colax $G$-morphisms} in which the sense of 
$\overline{f}$ is reversed. In the case of strict $G$-coalgebras,
this could be obtained from lax $G$-morphisms via a formal
dualization process, but in our more general setting this is
no longer the case. We shall write $\tilde{f}$ for the 2-cell
part of a colax $G$-morphism; the coherence conditions are 
\begin{align*}
  \vcenter{
\xymatrix{
& GA \ar[r]^{Ga} \ar[dr]_{Gf} \ar@{=>}[dd]^{\tilde{f}} & 
G^2A \ar[dr]^{G^2f} \ar@{=>}[d]^{G\tilde{f}} \\
A \ar[ur]^{a} \ar[dr]_{f} && GB \ar[r]_{Gb} \ar@{=>}[d]^{\beta} & G^2B \\
& B \ar[r]_{b} \ar[ur]_{b} & GB \ar[ur]_{dB} 
}}
&=
  \vcenter{
\xymatrix{
& GA \ar[r]^{Ga} \ar@{=>}[d]^{\alpha} & G^2A \ar[dr]^{G^2f}  \\
A \ar[ur]^{a} \ar[dr]_{f} \ar[r]^{a} & 
GA \ar[ur]_{dA} \ar[dr]^{Gf} \ar@{=>}[d]^{\tilde{f}} && G^2B \\
& B \ar[r]_{b} & GB \ar[ur]_{dB} 
}}
\\
  \vcenter{
\xymatrix{
&  & A \ar[dr]^{f}  \\
A  \ar[dr]_{f} \ar@/^1pc/[urr]^{1}_{~}="1" & 
&& B \\
& B \ar[r]_{b} \ar@/^1pc/[urr]^{1}_(0.53){~}="2" \ar@{=>}"2";[r]_{\beta_0} & 
GB \ar[ur]_{eB} 
}}
&=
  \vcenter{
\xymatrix{
&  & A \ar[dr]^{f}  \\
A  \ar[dr]_{f} \ar[r]^{a} \ar@/^1pc/[urr]^{1}_{~}="1" \ar@{=>}"1";[r]^{\alpha_0} & 
GA \ar[ur]_{eA} \ar[dr]^{Gf} \ar@{=>}[d]^{\tilde{f}} && B \\
& B \ar[r]_{b} & GB \ar[ur]_{eB} 
}}
\end{align*}

A $G$-transformation between 
lax $G$-morphisms $(f,\overline{f})$ and $(g,\overline{g})$ from 
$(A,a,\alpha,\alpha_0)$ to $(B,b,\beta,\beta_0)$ is a 2-cell
$\rho\cl f\to g$ satisyfying the single condition
\begin{align*}
  \vcenter{
\xymatrix{
A \ar@/^1pc/[r]^{f}_{~}="1" \ar[d]_{a} & B \ar[d]^{b} \\
GA \ar@/^1pc/[r]^{Gf}="2"_{~}="3" \ar@/_1pc/[r]_{Gg}^{~}="4" & GB
\ar@{=>}"1";"2"^{\overline{f}}
\ar@{=>}"3";"4"^{G\rho} 
}}
&=
  \vcenter{
\xymatrix{
A \ar@/^1pc/[r]^{f}_{~}="1" \ar@/_1pc/[r]^{~}="2"_{g}="3" \ar[d]_{a} & B \ar[d]^{b} \\
GA  \ar@/_1pc/[r]_{Gg}^{~}="4" & GB
\ar@{=>}"3";"4"^{\overline{g}}
\ar@{=>}"1";"2"^{\rho} 
}}
\end{align*}

There is a 2-category \LaxG of lax $G$-coalgebras, lax
$G$-morphisms, and $G$-transformations. 

There are also $G$-transformations between colax $G$-morphisms,
and a corresponding 2-category,
but we shall not need these.

\begin{example}
If $G$ is the 2-comonad on $[\ob X,\Cat]$ whose Eilenberg-Moore
2-category is $[X,\Cat]$, then a lax $G$-coalgebra is a lax functor
from $X$ to \Cat,
a lax $G$-morphism is a lax natural transformation, and a 
colax $G$-morphism is an oplax natural transformation. 
A $G$-transformation is a modification between oplax natural
transformations.  In particular, if $G$ is the identity 2-comonad
on \Cat, a lax $G$-coalgebra is a monad, lax $G$-morphism is 
a monad functor, a colax $G$-morphism is a monad opfunctor,
and a $G$-transformation is a monad transformation.
\end{example}

Generalizing the fact that every adjunction induces a monad
we have:

\begin{theorem}
  Let $(B,b)$ be a strict $G$-coalgebra, and let  
$u\cl B\to A$ be a morphism in \K with left adjoint $f\dashv u$,
and write $\eta$ and $\epsilon$ for the unit and counit. 
Then $A$ becomes a lax $G$-coalgebra $(A,a,\alpha,\alpha_0)$
where $a$, $\alpha$, and $\alpha_0$ are given by 
$$\xymatrix{
A \ar[r]^{f} & B \ar[r]^{b} & GB \ar[r]^{Gu} & GA }$$
$$\xymatrix{
A \ar[r]^{f} \ar[d]_{f} & B \ar[r]^{b} \ar[dl]^{1} & 
GB \ar[r]^{Gu}_(0.75){~}="1" \ar[dr]_{1}^(0.65){~}="2" & GA \ar[d]^{Gf} &
A \ar[d]_{f}  \ar@/^2pc/[dddrr]^{1}_(0.7){~}="7" \\
B \ar[d]_{b} \ar@{}[rrr]_{~}="5"  && {} & GB \ar[d]^{Gb} & B \ar[d]_{b} \ar[dr]^{1}_(0.37){~}="3" &  \\
GB \ar[d]_{Gu} \ar[rrr]_{dB}^{~}="6" && {} & G^2B \ar[d]^{G^2u} & 
GB \ar[d]_{Gu}  \ar[r]_{eB}^(0.3){~}="4" & B \ar[dr]_{u}^{~}="8" \\
GA \ar[rrr]_{dA} &  &  & G^2A & GA \ar[rr]_{eA} && A 
\ar@{=>}"1";"2"_{G\epsilon} 
\ar@{=>}"7";"8"_{\eta} 
}
$$
and then $u$ becomes a  lax $G$-morphism
$(u,\overline{u})\cl (B,b,\beta,\beta_0)\to(A,a,\alpha,\alpha_0)$
when we define $\overline{u}$ to be the 2-cell
$$\xymatrix{
B \ar[rrr]^{b} \ar[d]_{u} \ar[dr]^{1}_(0.4){~}="1" &&& GB \ar[d]^{Gu} \\
A \ar[r]_{f}^(0.3){~}="2" & B \ar[r]_{b} & GB \ar[r]_{Gu} & GA 
\ar@{=>}"2";"1"_{\epsilon} 
}$$
\end{theorem}

\proof
All of the assertions can be verified by direct calculation.
Alternatively, one can deduce the results from \cite{transport},
specifically Theorem~3.5 and the discussion preceding 
Theorem~4.1. As stated these results require 
the existence of comma objects in \K, but this can be avoided
by embedding \K in a larger 2-category if necessary.
\endproof

\begin{remark}
  This generalizes easily to the case where $B$ is only a 
lax $G$-coalgebra; see \cite{transport} again.
\end{remark}

By doctrinal adjunction (see \cite[Theorem~1.2]{Kelly-doctrinal})
the lax $G$-morphism structure on $u$ determines colax 
$G$-morphism structure on $f$ in the form of  the 2-cell
$\tilde{f}$ displayed below.
$$\xymatrix{
A \ar[r]^{f} \ar[d]_{f} & B \ar[r]^{b} & 
GB \ar[r]^{Gu}_(0.7){~}="1" \ar[dr]_{1}^(0.6){~}="2" & GA \ar[d]^{Gf} \\
B \ar[rrr]_{b} &&& GB 
\ar@{=>}"1";"2"_{G\epsilon} 
}$$

The Eilenberg-Moore construction in this context becomes the 
following
theorem. The lax descent objects mentioned in the theorem
were defined in \cite{codescent} via a minor modification 
of a definition in \cite{FibBicCorr}. They 
can be constructed out of inserters and equifiers.

\begin{theorem}\label{thm:classifier}
If $\K^G$ has lax descent objects for lax coherence data,
and so in particular if it has inserters and equifiers, then the
inclusion $\K^G\to\LaxG$ has a right adjoint.  
\end{theorem}

\proof 
This is a straightforward modification of the argument given
in \cite{codescent} for the case of lax $T$-algebras and lax
$T$-morphisms for a 2-monad $T$.
The value at a lax $G$-coalgebra $(A,a,\alpha,\alpha_0)$ of the
right adjoint is the universal strict $G$-coalgebra $\alg(A)$ 
equipped with a strict $G$-morphism $v\cl \alg(A)\to GA$
and a $G$-transformation 
$$\xymatrix{
\alg(A) \ar[r]^{v}_{~}="1" \ar[d]_{v} & GA \ar[d]^{Ga} \\
GA \ar[r]_{dA}^{~}="2" & G^2A 
\ar@{=>}"1";"2"^{\psi} }$$
satisfying the two conditions 
\begin{align}
\label{eq:descent1}
\vcenter{
\xymatrix{
& GA \ar[r]^{Ga} \ar[dr]^{Ga} \ar@{=>}[dd]^{\psi} & 
G^2A \ar[dr]^{G^2a} \ar@{=>}[d]^{G\alpha} \\
\alg(A) \ar[ur]^{v} \ar[dr]_{v}   && G^2A \ar[r]^{GdA} & G^3A \\
& GA \ar[ur]^{dA} \ar[r]_{dA} & G^2A \ar[ur]_{dGA} }}
&=
\vcenter{
\xymatrix{
& GA \ar[r]^{Ga} \ar@{=>}[d]^{\psi} & G^2A \ar[dr]^{G^2a} \\
\alg(A) \ar[ur]^{v} \ar[dr]_{v} \ar[r]^{v} & 
GA \ar[ur]^{dA} \ar[dr]^{Ga} \ar@{=>}[d]^{\psi}  && G^3A \\
& GA \ar[r]_{dA} & G^2A \ar[ur]_{dGA} }}
\\
\label{eq:descent2}
\vcenter{
\xymatrix{
& GA \ar[dr]_{Ga} \ar@{=>}[dd]^{\psi} \ar@/^1pc/[drr]^{1}_(0.49){~}="1" \ar@{=>}"1";[dr]^{G\alpha_0} \\
\alg(A) \ar[ur]^{v} \ar[dr]_{v} & & G^2A \ar[r]_{GeA} & GA \\
& GA \ar[ur]^{dA} \ar@/_1pc/[urr]_{1} }}
&= id
\end{align}
The component at $(A,a,\alpha,\alpha_0)$ of the counit is the 
lax morphism $(q,\overline{q})\cl \alg(A)\to(A,a,\alpha,\alpha_0)$
where $q$ is the composite $eA.v$ and $\overline{q}$ is given by
the pasting composite
$$\xymatrix{
\alg(A) \ar[dd]_{a'} \ar[dr]_{v} \ar[r]^{v} & 
GA \ar[rr]^{eA} \ar[dr]^{Ga} \ar@{=>}[d]^{\psi} && A \ar[dd]^{a} \\
& GA \ar[r]^{dA} \ar[dr]_{dA} & G^2A \ar[dr]_{eGA} \\
G\alg(A) \ar[rr]_{Gv} & & G^2A \ar[r]_{GeA} & GA }$$
in which $a'$ denotes the coalgebra structure of $\alg(A)$.  
\endproof

\begin{theorem}\label{thm:classifier-adjunction}
If \K has lax descent objects for lax coherence data, and
$G$ preserves them, then the inclusion $\K^G\to\LaxG$ has
a right adjoint. Furthermore, the  component at a lax
coalgebra $(A,a,\alpha,\alpha_0)$ of the counit of the adjunction
is a lax morphism
$(q,\overline{q})\cl (\alg(A),a')\to (A,a,\alpha,\alpha_0)$
for which $q\cl \alg(A)\to A$ has a left adjoint $p$, and this
$p$ inherits a colax $G$-morphism structure. 
\end{theorem}

\proof
The first sentence holds because the forgetful 2-functor
$U\cl \K^G\to\K$ creates limits of any type which are 
preserved by $G$. Thus we may construct $\alg(A)$ as the 
universal object in \K equipped with a morphism $v\cl \alg(A)\to GA$
and 2-cell $\psi$ as in the proof of the previous theorem. 
There is then a unique induced morphism $a'\cl \alg(A)\to G\alg(A)$
making $\alg(A)$ into a strict $G$-coalgebra, $v$ into a 
strict morphism of $G$-coalgebras, and $\psi$ into a 
$G$-transformation, and the construction of $q$ and
$\overline{q}$ was described in the previous theorem. 

The existence of the adjoint $p$ is a lax version of 
 \cite[Theorem~3.2]{codescent}, while the colax structure
of $p$ is given by doctrinal adjunction again: see 
\cite[Theorem~1.2]{Kelly-doctrinal}. Explicitly, $p\cl A\to\alg(A)$ 
is the unique morphism in \K whose composite with $v$
is $a\cl A\to GA$ and whose composite with $\psi$ is $\alpha$. 
The unit 
$\eta\cl 1\to qp$ of the adjunction $p\dashv q$ is $\alpha_0$. 
The counit is the unique 2-cell $\epsilon\cl qp\to 1$ whose 
composite $v\epsilon$ with $v$ is given by 
$$\xymatrix{
vpq \ar@{=}[r] & aq \ar@{=}[r] & a.eA.v \ar@{=}[r] & eGA.Ga.v \ar[d]^{eGa.\psi} \\
v \ar@{=}[rrr] &&& eGA.dA.v }$$
The colax structure for $p$ is given by the unique 2-cell
$\tilde{p}$ for which the pasting composite 
$$\xymatrix{
A \ar[r]^{p}_{~}="1" \ar[d]_{a} \ar@/^2pc/[rr]^{a} &
\alg(A) \ar[r]^{v} \ar[d]^{a'} & GA \ar[d]^{dA} \\
GA \ar[r]_{Gp}^{~}="2" \ar@/_2pc/[rr]_{Ga} & G\alg(A) \ar[r]_{Gv} &
G^2A 
\ar@{=>}"2";"1"^{\tilde{p}} }$$
is equal to $\alpha$.
\endproof

Under the hypotheses of the theorem, suppose that $(B,b)$ is 
a strict $G$-coalgebra, and that $u\cl B\to A$ is a morphism with a
left adjoint $f\dashv u$. Let $(A,a,\alpha,\alpha_0)$ be 
the induced lax $G$-coalgebra, and $(u,\overline{u})\cl (B,b)\to(A,a,\alpha,\alpha_0)$ the induced lax $G$-morphism. By the universal
property of $\alg(A)$, this corresponds to a unique strict 
$G$-morphism $k\cl (B,b)\to(\alg(A),a')$, with 
$(q,\overline{q})k=(u,\overline{u})$. 

\begin{definition}
We call this $k$ the
{\em canonical comparison}, and say that the original adjunction
$f\dashv u$, or just $u$, is {\em strictly monadic relative to $G$} if this $k$
is an isomorphism, and {\em monadic relative to $G$} if 
the underlying morphism $k\cl B\to\alg(A)$ in \K is an equivalence. 
\end{definition}

We now generalize the centipedes of \cite{Street-twoconstructions}
to this setting; of course they have long since lost any similarity
to centipedes in the biological sense. The description is a sort 
of dual of that in Section~\ref{sec:pic} since our forgetful 
2-functor is now comonadic rather than monadic. 

For a 2-category \M and an object $C\in\M$ we write $C\ls\M$ 
for the evident 2-category whose objects are arrows $A\cl C\to \partial_1 A$ in \M, whose arrows from $(A,\partial_1A)$ to $(B,\partial_1B)$ are 
arrows $\partial_1f\cl \partial_1 B\to \partial_1 A$ equipped with a 2-cell 
$f\cl A\to \partial_1f.B$,  as in 
\begin{align*}
  \vcenter{
\xymatrix{
& \partial_1A \\
C \ar[r]_{B}^(0.7){~}="2" \ar[ur]^{A}_(0.6){~}="1"  & \partial_1B \ar[u]_{\partial_1f} 
\ar@{=>}"1";"2"_{f} }}
\end{align*}
and whose 2-cells from $(f,\partial_1f)$ to $(g,\partial_1g)$  
are 2-cells from $\partial_1f\to \partial_1g$ satisfying the 
evident compatibility condition. There is an evident 
projection $\partial_1\cl(\C\ls\M)\to\M\op$.

We shall be interested in the case where $\M=\K^G$.

For $G$-coalgebras $(C,c)$ and $(B,b)$, we define a
{\em $(C,c)$-centred centipede in $(B,b)$} to be a parallel pair 
of arrows in $(C,c)\ls\K^G$ which lie over $Gb\cl GB\to G^2B$
and $dB\cl GB\to G^2B$. We display below such a pair in two
ways.

\begin{align*}
\vcenter{
\xymatrix @C6pc @R3pc {
& G^2B \\
C \ar[ur]^{M}_(0.5){~}="1"_(0.7){~}="2" \ar[r]_{N}^(0.53){~}="3"^(0.73){~}="4" & GB \ar@<1ex>[u]^{Gb} \ar@<-1ex>[u]_{dB} 
\ar@{=>}"1";"3"^{m} \ar@{=>}"2";"4"^{n} 
}}
\quad
\vcenter{
\xymatrix{
(M,G^2B) \ar@<1ex>[r]^{(m,Gb)} \ar@<-1ex>[r]_{(n,dB)} & (N,GB) 
}}
\end{align*}
Rather than {\em universal reflection}, we shall simply speak
of a {\em coequalizer} of a centipede, and define this to be  
a coequalizer
in $(C,c)\ls\K^G$ which is 
preserved by $\partial_1$
(and so must be projected to $b\cl B\to GB$). 

A morphism $g\cl(B,b)\to(B',b')$ in $\K^G$ sends a 
centipede in $(B,b)$ to one in $(B',b')$, but it need not 
of course preserve coequalizers. A coequalizer
of a centipede in $(B,b)$ 
is said to be {\em absolute} when it is preserved by all morphisms
$g\cl(B,b)\to(B',b')$ in $\K^G$.

Among the absolute coequalizers are the split ones, and we
can similarly define split coequalizers of centipedes.
But this would require that  the diagram
$$\xymatrix{
B \ar[r]^{b} & GB \ar@<1ex>[r]^{Gb} \ar@<-1ex>[r]_{dB} & G^2B}$$
which is always an equalizer in $\K^G$, is in fact a split
equalizer in $\K^G$ (not just in \K). This is not true for a general
coalgebra $(B,b)$, but it is true for a cofree one, such as
$(GA,dA)$, so we shall only consider split centipedes in 
cofree coalgebras. We always use the standard splitting for
cofree coalgebras, as in the following diagram.
$$\xymatrix{
GA \ar[r]^{dA} & 
G^2A \ar@<1ex>[r]^{GdA} \ar@<-1ex>[r]_{dGA} \ar@/^1pc/[l]^{GeA} & 
G^3A \ar@/^2pc/[l]^{G^2eA} }$$

Explicitly, then, a split centipede in $(GA,dA)$ has the 
form 
$$\xymatrix{
(M,G^3A) \ar@<1ex>[r]^{(m,GdA)} \ar@<-1ex>[r]_{(n,dGA)} &
(N,G^2A) \ar[r]^{(p,dA)} \ar@/^2pc/[l]^{(t,G^2eA)} &
(P,GA) \ar@/^1pc/[l]^{(s,GeA)} }$$ 
where 
$(p,dA)(m,GdA) = (p,dA)(n,dGA)$, 
$(p,dA)(s,GeA) =1$, 
$(s,GeA)(p,dA) = (n,dGA)(t,G^2eA)$, and 
$(m,GdA)(t,G^2eA) =1$.
Clearly every split coequalizer of a centipede is 
absolute. 

For a morphism $g:(B,b)\to(B',b')$ in $\K^G$, we say that a coequalizer of a
centipede in $(B,b)$ is {\em $g$-split} or {\em $g$-absolute}
if the induced centipede in $(B',b')$ is so.

We are now ready to state our version of the Beck monadicity
theorem.

\begin{theorem}\label{thm:strict-monadicity}
Suppose that  \K is a 2-category admitting lax descent objects of coherence data and that $G=(G,d,e)$ is a 2-comonad on \K for which 
the 2-functor $G$ preserves such lax descent objects; for example,
\K could have inserters and equifiers, and $G$ preserve them. 
Let  $(B,b)$ be a strict $G$-coalgebra, and $u\cl B\to A$ a 
right adjoint in \K, and let $w\cl B\to GA$ be the unique strict
$G$-morphism whose composite with $eA\cl GA\to A$ is $u$.
The following conditions are equivalent:
\begin{enumerate}[(a)]
\item $u$ is strictly monadic relative to $G$
\item $w$ creates $w$-absolute coequalizers of centipedes
\item  $w$ creates coequalizers of $w$-split centipedes in $(B,b)$.
\end{enumerate}
\end{theorem}

\proof
We shall prove $(a)\Rightarrow(b)\Rightarrow(c)\Rightarrow(a)$,
following the classical proof as modified by Par\'e to use absolute
coequalizers. Of course the implication $(b)\Rightarrow(c)$ holds
trivially since split coequalizers are absolute. 

First suppose (a). Suppose that the pair $(m,Gb)\cl (M,G^2B)\to(N,GB)$ and $(n,dB)\cl (M,G^2B)\to(N,GB)$ are sent by $w$ to a 
centipede in $(GA,dA)$ with an absolute coequalizer
$$\xymatrix{
(G^2w.M,G^3A) \ar@<1ex>[r]^{(Gw.m,GdA)} \ar@<-1ex>[r]_{(Gw.n,dGA)} & (Gw.N,G^2A) \ar[r]^-{(q,dA)} & (Q,GA).
}$$
Then $Ga$ preserves the coequalizer, and so  
the rows of the diagram
$$\xymatrix @C3pc {
(G^3a.G^2w.M,G^4A) \ar@<1ex>[r]^{(G^3a.Gw.m,GdGA)} \ar@<-1ex>[r]_{(G^3a.Gw.n,dG^2A)} \ar[d]_{(G^2\psi.M,G^2dA)} & 
(G^2a.Gw.N,G^3A) \ar[r]^-{(G^2a.q,dGA)} \ar[d]_{(G\psi.N,GdA)} & 
(Ga.Q,G^2A) \ar@{.>}[d]^{(\theta,dA)} \\
(G^2w.M,G^3A) \ar@<1ex>[r]^{(Gw.m,GdA)} \ar@<-1ex>[r]_{(Gw.n,dGA)} &
(Gw.N,G^2A) \ar[r]_{(q,dA)} & (Q,GA)
}$$
are  coequalizers. It follows by the universal property of the
top coequalizer that there is a unique $\theta$ as in the dotted
arrow which makes the whole diagram commute. By
considering the coequalizer obtained by applying $G^2a$ to the
top row, one verifies that this $\theta$ satisfies the coherence
condition \eqref{eq:descent1}, and by using the original 
coequalizer one checks \eqref{eq:descent2}; thus there is a 
unique map $P\cl C\to B$ in $\K^G$ with $wP=Q$ and 
$\psi P=q$.

Commutativity of the square defining $\theta$ is exactly 
what is needed to apply the two-dimensional aspect of the
universal property of $G\alg(A)$ (recalling also that $G$ preserves
the limit which defines $\alg(A)$), and so there is a unique 
induced 2-cell $p\cl N\to P$ whose composite with $w\cl B\to GA$
is $q$. One now follows the same steps as in the classical case
to show that 
$$\xymatrix{
(M,G^2B) \ar@<1ex>[r]^{(m,Gb)} \ar@<-1ex>[r]_{(n,dB)} & (N,GB) \ar[r]^-{(p,b)} & (P,B)
}$$
is the required coequalizer. This proves (b).

Thus it remains to show that (c) implies (a). As in the classical 
proof, we show 
that the comparison $k\cl B\to \alg(A)$ is invertible by 
constructing a suitable centipede with a $w$-split coequalizer.  
The centipede is centred at $\alg(A)$, and is displayed below.
$$\xymatrix{
(G^2f.Ga.v,G^2B) \ar@<1ex>[r]^{(G\tilde{f}.v,Gb)}
\ar@<-1ex>[r]_{(G^2f.\psi,dB)} & (Gf.v,GB) }$$
Applying the strict map $w\cl B\to GA$, and using the fact
that $wf=a$ and $Gw.\tilde{f}=\alpha$, we see that the induced
centipede is given as in the following diagram
$$\xymatrix{
(G^2a.Ga.v,G^3A) \ar@<1ex>[r]^{(G\alpha.v,GdA)}
\ar@<-1ex>[r]_{(G^2a.\psi,dGA)} & 
(Ga.v,GdA) \ar[r]^-{(\psi,dA)} \ar@/^2pc/[l]^{(G^2\alpha_0.Ga.v,G^2eA)} & 
(v,dA) \ar@/^1pc/[l]^{(G\alpha_0.v,GeA)} 
}$$
in which a split coequalizer for the induced centipded is also given. The lifted coequalizer now has the form 
$$\xymatrix{
(G^2f.Ga.v,G^2B) \ar@<1ex>[r]^{(G\tilde{f}.v,Gb)}
\ar@<-1ex>[r]_{(G^2f.\psi,dB)} & (Gf.v,GB) \ar[r]^-{(p,b)}
& (P,B)
}$$
where $P\cl\alg(A)\to B$  satisfies $wP=v$ and 
$G\overline{u}.b.P=\psi$. One now uses the universal 
property of $\alg(A)$ to check that $kP=1$, and the 
uniqueness aspect in the creation of $w$-absolute coequalizers
to show that $Pk=1$. 
\endproof

We leave to the reader the modifications necessary to prove
the following ``up-to-equivalence'' version of monadicity.

\begin{theorem}\label{thm:monadicity}
In the setting of Theorem~\ref{thm:strict-monadicity} once again,
the following conditions are equivalent:
\begin{enumerate}[(a)]
\item $u$ is monadic relative to $G$
\item any centipede in $(B,b)$ which is sent by $w$ to a
centipede in $(GA,dA)$ with an absolute coequalizer itself has 
a coequalizer, and such coequalizers are preserved and reflected
by $w$
\item as in (b) but with split coequalizers rather than absolute
ones.
\end{enumerate}
\end{theorem}

\section{Dualization}

One of the striking things about the paper \cite{ftm} was the 
use of duality to obtain both Eilenberg-Moore objects and Kleisli 
objects for both monads and comonads, all in a single setting.

In this short final section, we very briefly indicate how the 
results of the previous section can be dualized. 

Write $\K\op$ for the 2-category obtained from \K by formally
reversing the 1-cells, but not the 2-cells. A 2-comonad $G$
on \K, as in the previous section, corresponds to a 2-monad $T$
on $\K\op$. A lax $G$-coalgebra is then the same thing as a 
lax $T$-{\em{algebra}}. Given lax $G$-coalgebras $A$ and $B$, a 
lax morphism of $G$-coalgebras  from $A$ to $B$ is the same thing
as a {\em colax} morphism $T$-algebras from $B$ to $A$. 
Thus $\LaxG=\LaxT\op$.

The assumption of certain limits in \K, preserved by $G$, is 
equivalent to the existence of certain {\em colimits} in $\K\op$,
preserved by $T$. 

A right adjoint to the inclusion $\K^G\to\LaxG$ is equivalent to 
a {\em left} adjoint to the inclusion $\K^T\to\LaxT$ of the 
Eilenberg-Moore 2-category $\K^T$ of $T$ into the 2-category
$\LaxT$ of lax $T$-algebras, colax $T$-morphisms, and 
$T$-transformations. 

Thus our Beck-style theorem of the previous section becomes
a recognition theorem for colax morphism classifiers (in the 
context of lax algebras for a 2-monad). 

If the original comonad $G$ is the identity, so that $T$ is 
also the identity, then these colax morphism classifiers
are in fact Kleisli objects for monads. 

Write $\K\co$ for the 2-category obtained from \K by formally
reversing the 2-cells but not the 1-cells. A 2-comonad $G$ on \K
can equally be seen as a 2-comonad on $\K\co$, but now lax
$G$-coalgebras in \K are the same as {\em colax} $G$-coalgebras
in $\K\co$.

If $G$ (and so $T$) are identities, then this reduces to 
Eilenberg-Moore objects for {\em comonads}.

One can also reverse both 1-cells and 2-cells to get a 
2-category $\K\coop$, and the theory now generalizes 
Kleisli objects for comonads.

\bibliographystyle{plain}

\begin{thebibliography}{10}

\bibitem{bicategories}
Jean B{\'e}nabou.
\newblock Introduction to bicategories.
\newblock In {\em Reports of the Midwest Category Seminar}, pages 1--77.
  Springer, Berlin, 1967.

\bibitem{BKP}
R.~Blackwell, G.~M. Kelly, and A.~J. Power.
\newblock Two-dimensional monad theory.
\newblock {\em J. Pure Appl. Algebra}, 59(1):1--41, 1989.

\bibitem{BohmMenini-pretorsors}
Gabriella B{\"o}hm and Claudia Menini.
\newblock Pre-torsors and {G}alois comodules over mixed distributive laws.
\newblock {\em Appl. Categ. Structures}, 19(3):597--632, 2011.

\bibitem{BrzezinskiMarquezVercruysse}
Tomasz Brzezi{\'n}ski, Adrian Vazquez~Marquez, and Joost Vercruysse.
\newblock The {E}ilenberg-{M}oore category and a {B}eck-type theorem for a
  {M}orita context.
\newblock {\em Appl. Categ. Structures}, 19(5):821--858, 2011.

\bibitem{ElKaoutit}
L.~El~Kaoutit.
\newblock Wide {M}orita contexts in bicategories.
\newblock {\em Arab. J. Sci. Eng. Sect. C Theme Issues}, 33(2):153--173, 2008.

\bibitem{Kelly-doctrinal}
G.~M. Kelly.
\newblock Doctrinal adjunction.
\newblock In {\em Category Seminar (Proc. Sem., Sydney, 1972/1973)}, pages
  257--280. Lecture Notes in Math., Vol. 420. Springer, Berlin, 1974.

\bibitem{transport}
G.~M. Kelly and Stephen Lack.
\newblock Monoidal functors generated by adjunctions, with applications to
  transport of structure.
\newblock In {\em Galois theory, {H}opf algebras, and semiabelian categories},
  volume~43 of {\em Fields Inst. Commun.}, pages 319--340. Amer. Math. Soc.,
  Providence, RI, 2004.

\bibitem{codescent}
Stephen Lack.
\newblock Codescent objects and coherence.
\newblock {\em J. Pure Appl. Algebra}, 175(1-3):223--241, 2002.

\bibitem{Pecsi}
Bertalan P{\'e}csi.
\newblock On {M}orita contexts in bicategories.
\newblock {\em Appl. Categ. Structures}, 20(4):415--432, 2012.

\bibitem{ftm}
Ross Street.
\newblock The formal theory of monads.
\newblock {\em J. Pure Appl. Algebra}, 2(2):149--168, 1972.

\bibitem{Street-twoconstructions}
Ross Street.
\newblock Two constructions on lax functors.
\newblock {\em Cahiers Topologie G\'eom. Diff\'erentielle}, 13:217--264, 1972.

\bibitem{FibBicCorr}
Ross Street.
\newblock Correction to: ``{F}ibrations in bicategories'' [{C}ahiers
  {T}opologie {G}\'eom.\ {D}iff\'erentielle {\bf 21} (1980), no.\ 2, 111--160;
  {MR}0574662 (81f:18028)].
\newblock {\em Cahiers Topologie G\'eom. Diff\'erentielle Cat\'eg.},
  28(1):53--56, 1987.

\end{thebibliography}

\end{document}